\documentstyle[11pt]{article}

\newcommand{\noi}{\noindent}
\newcommand{\de}{\backslash}

\newcommand{\qed}{\hfill $\Box $}
\newcommand{\pf}{\noindent {\bf Proof.} }

\vspace{-0.3in}

\begin{document}
\hbox{Discrete Math. 306(2006), no.\,15, 1798--1804.}
\bigskip

\centerline{\bf On the unique representability of spikes over prime fields}
\bigskip
\centerline{Zhaoyang Wu and Zhi-Wei Sun}
 \medskip
\centerline{Department of Mathematics, Nanjing University}
\centerline{Nanjing 210093, People's Republic of China}
\centerline{{\tt wu$_{-}$zy@nju.edu.cn\ \ zwsun@nju.edu.cn}}

\date{}

%\begin{document}
%\maketitle

\begin{abstract}  For an integer $n \ge 3$, a rank-$n$ matroid is called an
$n$-spike if it consists of $n$ three-point lines through a
common point such that, for all $k\in\{ 1, 2, \ldots, n - 1\}$,
the union of every set of $k$ of these lines has rank $k+1$.
Spikes are very special and important in matroid theory. Wu
\cite{wu03} found the exact numbers of $n$-spikes over fields with
2, 3, 4, 5, 7 elements, and the asymptotic values for larger
finite fields. In this paper, we prove that, for each prime number
$p$, a $GF(p$) representable $n$-spike $M$ is only representable on
fields with characteristic $p$ provided that $n \ge 2p-1$.
Moreover, $M$ is uniquely representable over $GF(p)$.

\vspace{0.2in}
 \noi {\it Keywords}: Matroid, Spike, Unique Representability

\vspace{0.2in}
\end{abstract}

\section{Introduction}

\hspace{.6cm} Spikes are special and important matroids. They are
appearing with increasing frequency in the
 matroid theory literature. Long before the name "spike"
was introduced, the Fano and non-Fano matroids, two examples of
$3$--spikes, had already appeared in almost every corner of
matroid theory \cite{ox92}. Oxley [6, Section 11.2] showed that
all rank--$n$, $3$--connected binary matroids without a $4$--wheel
minor can be obtained from a binary $n$--spike by deleting at most
two elements. Oxley, Vertigan, and Whittle \cite{ov96} used spikes
and one other class of matroids to show that, for all $q \ge 7$,
there is no fixed bound on the number of inequivalent
$GF(q)$--representations of a $3$--connected matroid, thereby
disproving a conjecture of Kahn \cite{ka88}.

\vspace{-0.2cm} Ding, Oporowski, Oxley, and Vertigan \cite{doov97}
showed that every sufficiently large $3$--connected matroid has,
as a minor, $U_{2, n + 2}$, $U_{n, n + 2}$, a wheel or whirl of
rank $n$, $M(K_{3,n})$, $M^*(K_{3,n})$, or an $n$--spike.
Moreover, Wu \cite{wu97} showed that spikes, like wheels and whirls,
can be characterized in terms of a natural extremal connectivity
condition. Wu \cite{wu03} discussed the representability of spikes
over finite fields, and found the exact numbers of $n$-spikes over
fields with at most seven elements, and the asymptotic values for
larger finite fields . One referee for the last mentioned paper
was interested in the problem that on what conditions a
$GF(p)$-representable spike is only representable over fields with
characteristic $p$. We consider this problem an interesting one
with fair importance in matroid theory, and this paper is a
response to the problem.

For $n\ge 3$, a matroid $M$ is called an {\it $n$--spike with tip
$t$} \cite{doov97} if it satisfies the following three conditions:
\begin{enumerate}{\it
\item[{\rm (i)}] The ground set is the union of $n$ lines, $L_1$, $L_2$, \ldots,
$L_n$, all having three points and passing through a common point
$t$;
\item[{\rm (ii)}]  \(r( \cup_{i=1}^{k}{L_i}) = k+1\) for all $k \in \{ 1,2, \ldots, n-1 \}$; and
\item[{\rm (iii)}] \(r({L_1} \cup {L_2} \cup \ldots \cup {L_n}) = n\)}.
\end{enumerate}

\noi In this paper, an {\it $n$--spike with tip $t$} will be
simply called an {\it $n$--spike}.

Some $3$--spikes have the property that more than one element may
be viewed as the tip of the spike. However, it is clear that the
tip is unique for an $n$--spike when $n \ge 4$. Since there are
only six $3$--spikes, and it is easy to verify all our results for
the case $n = 3$, we will assume that $n$ is at least four in the
proofs of our lemmas and theorems so that we can fix the tip.

For an $n$--spike $M$ representable over a field $F$, if we choose
a base $\{ 1, 2, \ldots, n \}$ containing exactly one element from
each of the lines $L_i$, then $M$ can be represented in the form
\vspace{0.14cm}

\hspace{1.5cm} $\begin{array} {l} 1 \\ 2 \\ 3 \\ \vdots \\ n
\end{array}
  \left[ \begin{array}{cccccccccccc}
   1 & 0 & 0 & \ldots & 0 & | & 1 & {1 + x_1} & 1 & 1 & \ldots & 1 \\
   0 & 1 & 0 & \ldots & 0 & | & 1 & 1 & {1 + x_2} & 1 & \ldots & 1 \\
   0 & 0 & 1 & \ldots & 0 & | & 1 & 1 & 1 & {1 + x_3} & \ldots & 1 \\
   \vdots & \vdots & \vdots & \ddots & \vdots & | & \vdots & \vdots & \vdots
& \vdots & \ddots & \vdots \\
   0 & 0 & 0 & \ldots & 1 & | & 1 & 1 & 1 & 1 & \ldots & {1 + x_n}
\end{array} \right],$

\vspace{0.1cm} \noi where the tip of $M$ corresponds to column $n
+ 1$. We shall call this matrix a {\it special standard
representation} of $M$ and $\{ 1, 2, \ldots, n \}$ the {\it
distinguished basis} associated with the representation. Clearly,
this matrix is uniquely determined by the vector $( x_1, x_2,
\ldots, x_n )$. We shall call this vector the {\it diagonal} of
the representation.

Two matrix representations  $A_1$ and $A_2$ are {\it equivalent}
if $A_1$ can be obtained from $A_2$ by a sequence of the following
six operations. (For details, see [5, Section 6.3].)

\begin{enumerate}{\it
\item[{\rm (i)}] Interchange two rows.
\item[{\rm (ii)}] Scale a row, that is, multiply it by a non-zero member of $F$.
\item[{\rm (iii)}] Replace a row by the sum of that row and another.
\item[{\rm (iv)}] Interchange two columns (moving their labels with the columns).
\item[{\rm (v)}] Scale a column, that is, multiply it by a non-zero member of $F$.
\item[{\rm (vi)}] Replace each entry of the matrix by its image under some
automorphism of $F$.}
\end{enumerate}

$A_1$ and $A_2$ are {\it weakly equivalent} if we are also allowed
to relabel the matroid, that is, $A_1$ can be obtained from $A_2$
by a sequence of operations {\it {\rm (i) - (vii)}} where the last
of these operations is the following:

\hspace{-0.5cm} {\it {\rm (vii)} Relabel the columns.}

Since our main purpose is to study the conditions on which a
matroid is or is not representable over a finite field, we will
often consider unlabeled matroids.
 Thus, we will frequently ignore the labels on elements of matroids, and consider
weak equivalence.

If two special standard representations are weakly equivalent,
their corresponding diagonals will also be said to be {\it weakly
equivalent}. Two diagonals are {\it distinct} if they are not
weakly equivalent. Two elements of an $n$--spike are {\it
conjugate} if they lie on the same line $L_i$ and neither of them
is the tip. In a special standard representation of a given spike,
if we interchange some base elements with their conjugates, and
standardize the resulting matrix, we obtain another special
standard representation of the spike. Moreover, all possible
special standard representations of the spike are obtainable in
this way. In the rest of the paper, we shall call this
interchanging-standardizing procedure {\it swapping}. For two
special standard representations $A_1$ and $A_2$ of an $n$-spike,
the distinguished bases of $M[A_1]$ and $M[A_2]$ are $n$--element
subsets intersecting all the lines $L_i$. Since the tip is fixed
and is in neither distinguished basis, $A_1$ and $A_2$ are weakly
equivalent if and only if we can obtain the distinguished
 basis of $M[A_1]$ from that
of $M[A_2]$ by swappings. Therefore, $A_1$ and $A_2$ are weakly
equivalent if and only if $A_1$ can be obtained from $A_2$ by a
sequence of swappings, and replacing each entry of the resulting
matrix by its image under some automorphism of the field $F$.

In the rest of this paper, the matroid notation and terminology
will follow Oxley \cite{ox92}. The notation and terminology for
spikes will follow Wu \cite{wu03}, of which some results are
quoted and some techniques are inherited in this paper.

\noindent{\bf (1.1) Main Theorem.} {\it For each prime number $p$,
if the integer $n$ is greater than or equal to $2p-1$, then an
$n$-spike that is $GF(p)$-representable can only be represented
over fields with characteristic $p$. Moreover, $M$ is uniquely
representable over $GF(p)$. }

\section{Preliminaries}

\hspace{.5cm} In the following sections, we use the notation
$[m,n]$ to denote the set of consecutive integers from $m$ to $n$,
namely, $\{ m, m+1, m+2, \ldots, n \}$, for our convenience.

 \noi {\bf
(2.1) Lemma.} {\it Let $p$ be a prime integer, $n$ be an integer
satisfying the condition $n \ge p-1$, and $a_1, a_2, \ldots, a_n
\in GF(p)\de \{ 0\}$. Suppose that $k \in GF(p)\de \{ 0\}$. Then
there is a non-empty subset $I\subseteq[1,n]$ such that

\hspace{5.25cm} $\sum_{i \in I}{a_i} = k.$}

\noi  \pf Viewing $GF(p)$ as $\textbf{Z}_p =
\textbf{Z}/p\textbf{Z}$, we rewrite $a_i$'s as

\hspace{5.25cm} $a_i = m_i+p\textbf{Z}$,

\noi where $1\le m_i \le p-1$.

As the system $\{ i+\frac{p}{m_i}\textbf{Z} \}_{i=1}^{p-1}$ covers
$\{ 1, 2, \ldots, p-1\}$, but not all the integers, by Theorem $1$
of Sun \cite{sun2} or Corollary $5$ of Sun \cite{sun1}, we can not
have $|\{ \{\sum_{i \in I}\frac{m_i}{p} \} : I \subseteq [1,p-1]
\}| \le p-1$. Therefore, $\{\sum_{i \in I}a_i : I \subseteq
[1,p-1] \} = \textbf{Z}_p$, and the lemma follows. \qed

\noi {\bf (2.2) Lemma.} {\it Let $p$ be a prime integer, $n$ be an
integer satisfying the condition $n \ge p$, and $a_1, a_2, \ldots,
a_n \in GF(p)$. Then there is a non-empty subset $I \subseteq
[1,n]$ such that

\hspace{5.25cm} $\sum_{i \in I}{a_i} = 0$.}

\noi  \pf Again the lemma can be easily derived from Theorem $1$
of Sun \cite{sun2}, and details are thus omitted . \qed

The following proposition is not hard to prove by induction. We
shall omit the proof.

\noi {\bf (2.3) Proposition.} {\it Suppose that $n$ is a positive
integer and that $x_i \neq 0$ for all $i \in \{ 1, 2, \ldots, n
\}$. Then the determinant of the matrix
  \[ det \left( \begin{array}{ccccc}
    {1 + x_1} & 1 & 1 & \ldots & 1 \\
    1 & {1 + x_2} & 1 & \ldots & 1 \\
    1 & 1 & {1 + x_3} & \ldots & 1 \\
    \vdots & \vdots & \vdots & \ddots & \vdots \\
    1 & 1 & 1 & \ldots & {1 + x_n}
\end{array} \right)  =  [1 + \sum^n_{i = 1}{x_i}^{-1}]\cdot \prod^n_{i = 1}x_i.\] }

Suppose that $A$ is a special standard representation of an
$n$--spike $M$ over $GF(p)$ with its diagonal $\vec{x}=( x_1, x_2,
\ldots, x_n )$, and that $C$ is a circuit-hyperplane of $M$.
Suppose that $I$ is the subset of $[1,n]$ such that the elements
of $M$ corresponding to $\{ x_i | i \in I \}$ is the intersection
of $C$ and the set of elements of $M$ to  diagonal $\vec{x}$. Then
we deduce by (2.3) that

\noi {\bf (2.4)} \hspace{4cm} $\sum_{i \in I}{x_i}^{-1} = -1$.

Conversely, suppose that $I$ is a subset of $[1,n]$ such that
equality (2.4) holds. Then the elements of $M$ corresponding to
$\{ x_i | i \in I \}$ combined with the conjugates of the elements
corresponding to the remainder of $\vec{x}$ form a
circuit-hyperplane of $M$. Therefore, every circuit-hyperplane $C$
corresponds to an $I$ with $\{ x_i | i \in I \}$ satisfying (2.4),
and vise versa. We denote the sub-vector corresponding to $\{ x_i
| i \in I \}$ by $\vec{x}(I)$, and call
 $C$ the $\it circuit-hyperplane$ $corresponding$  $to$ $\vec{x}(I)$.

\noi {\bf (2.5) Proposition.} {\it Let $A$ be a special standard
representation of an $n$--spike $M$ over $GF(p)$ with diagonal
$\vec{x}=( x_1, x_2, \ldots, x_n )$. Suppose that $n \ge p-1$.
Then $\vec{x}$ is weakly equivalent to a diagonal whose first
element is $-1$.}

\noi  \pf Since $n \ge p-1$,  we derived from (2.1) that there is
a sub-set $I \subseteq [1,n]$ that satisfies (2.4). By weakly
equivalence we may assume that there is a positive integer $m \le
n$ such that $I = [1,m]$. Let $C$ be the circuit-hyperplane
corresponding to $\vec{x}(I)$. By swapping all but the first
element of $\vec{x}(I)$ with their conjugates, we obtain a new
special standard representation of $M$. Let the diagonal of this
new special standard representation be $\vec{y} = ( y_1, y_2,
\ldots, y_n )$. It is obvious that for all elements corresponding
to $\vec{y}$ only the one corresponding to $y_1$ is contained in
$C$. The desired result thus follows by (2.3). \qed

\section{Proof of The Main Theorem}

\hspace{.6cm}We first introduce the following two propositions.

\noi {\bf (3.1) Proposition.} {\it Let $p$ be an odd prime
integer, $n$ be an integer with $n \ge 2p-1$, and matroid $M$ be
an $n$-spike representable over $GF(p)$ and another finite field
$F$ with characteristic $q$. Suppose that $A_1$, $A_2$ are two
special standard representations of $M$ over $GF(p)$ and $F$, and
their diagonals are $\vec{x} = ( -1, x_2, \ldots, x_n )$, and
$\vec{y} = ( y_1, y_2, \ldots, y_n )$, respectively. Suppose that
$m \in \textbf{Z} \de \{ 0 \}$, and $|m| \le (p-1)/2$, and
  $I$ is a subset of $[2,n]$ such that

\hspace{4cm} $|I| \le p-1$, and $\sum_{i \in I}{x_i}^{-1} = m$.

\noi Then we have the equality

\hspace{5.25cm} $\sum_{i \in I}{y_i}^{-1} = m$.}

\noi \pf Suppose that $C$ is the circuit-hyperplane corresponding
to $x_1 = -1$ of $\vec{x}$. Since $A_1$ and $A_2$ represent the
same spike, we deduce by (2.4) that $y_1$ of $\vec{y}$ is also
equal to $-1$.

Consider the case that $m = -1$. In this case, since $\sum_{i \in
I}{x_i}^{-1} = -1$, we consider the circuit-hyperplane
corresponding to $\vec{x}(I)$. For the reason that $A_1$ and $A_2$
are both special standard representation of the same spike, we
conclude that $\sum_{i \in I}{y_i}^{-1} = -1$.

Now consider the case that $m = 1$. Let $K = [2,n]\de I$. Since
$K$ has at least $p-1$ elements, we deduce by Lemma (2.1) that
there is a subset $L$ of $K$, such that

\hspace{5.25cm} $\sum_{i \in L}{x_i}^{-1} = -1$.

\noi Applying a discussion the same as that of the last paragraph,
we conclude that

\hspace{5.25cm} $\sum_{i \in L}{y_i}^{-1} = -1$.

\noi Let $I^{'} = I \bigcup L \bigcup \{ 1 \}$. Then we have

\hspace{5.25cm} $\sum_{i \in I^{'}}{x_i}^{-1} = -1$.

Therefore, there is a circuit-hyperplane $C$ corresponding to
$\vec{x}(I^{'})$. Since $A_1$ and $A_2$ are both special standard
representation of the same spike, we conclude that $\sum_{i \in
I^{'}}{y_i}^{-1} = -1$. It follows that $\sum_{i \in I}{y_i}^{-1}
= 1$.

Using the above result and the same technique, we can now prove
Proposition (3.1) for the case $m = -2$. Moreover, it is now clear
that we can complete the proof by induction. The details are thus
omitted. \qed

\noi {\bf (3.2) Lemma.} {\it Let $p$ be an odd prime integer, and
matroid $M$ be an $n$-spike representable over $GF(p)$. Suppose
that $n \ge 2p-1$. Then $M$ is uniquely representable over
$GF(p)$.}

\noi \pf Suppose that $A_1$, $A_2$ are two special standard
representations of $M$ over $GF(p)$, and their diagonals are
$\vec{x} = ( x_1, x_2, \ldots, x_n )$, and $\vec{y} = ( y_1, y_2,
\ldots, y_n )$.

We may assume, by (2.5), that $x_1 = -1$. For singleton set $I =
\{i\}$ with $x_i = m$, $m \in \textbf{Z} \de \{ 0 \}$ with $|m|
\le (p-1)/2$, we deduce by Proposition (3.1) that $y_i = x_i$.
Lemma (3.2) follows immediately. \qed

\noi {\bf Proof of the Main Theorem}  Since it is well known that
binary spikes are uniquely representable only on fields of
characteristic $2$, we only need to prove the main theorem with
odd prime number $p$.

Having Lemma (3.2) in hand, we only need to prove that $M$ is not
representable over field with characteristic not equal to $p$.
Suppose that $F$ is a field with characteristic $q$, and $M$ is
representable over $F$. We prove in the following that the prime
$q$ must be equal to $p$.

 Suppose that $A_1$, $A_2$ are special standard representations of $M$ over
$GF(p)$ and $F$, and $\vec{x} = ( x_1, x_2, \ldots, x_n )$, and
$\vec{y} = ( y_1, y_2, \ldots, y_n )$ are the diagonals
corresponding to $A_1$ and $A_2$, respectively. We assume, as we
may, that $x_1 = -1$. Moreover, we use values in $[-\frac{p-1}{2},
\frac{p-1}{2}] \de \{ 0 \}$ to represent the value of each
$x_i^{-1}$ of $\vec{x}$. Applying Proposition (3.1), we conclude
that, for each $i \in [1,n]$, $y_i^{-1} = x_i^{-1}$ in $GF(q)$.
Consequently, $M$ is representable over $GF(q)$, and we may assume
that $F = GF(q)$. As a result of the last assumption, we may
assume that $q \ge p$ in the following discussion.

 Now, consider the subscription set $I = [2,n]$. We partition $I$ into
 two parts $I_+$ and $I_-$, where $I_+ = \{ i \in [2,n] : x_i^{-1}> 0 \}$,
  and $I_- = \{ i \in [2,n] : x_i^{-1}< 0 \}$.

 First consider the case that $|I_-| \ge p$. We deduce, by Lemma (2.2),
 that there is a non-empty subset $L$ of $I_-$, such that

\hspace{5.25cm} $\sum_{i \in L}{x_i}^{-1} = 0$.

\noi This equality implies that $M$ has a circuit-hyperplane
corresponding to $\vec{x}(L\bigcup\{1\})$. Since $A_2$ is also a
special standard representation of $M$, we conclude that the
equality

\hspace{5.25cm} $\sum_{i \in L}{y_i}^{-1} = 0$

\noi holds in $GF(q)$. That is, the equality

\hspace{5.25cm} $\sum_{i \in L}{x_i}^{-1} = 0$

\noi holds in both $GF(p)$ and $GF(q)$.

Consider the sum $s = \sum_{i \in L}{x_i}^{-1}$ in $\textbf{Z}$.
 Since all values of $x_i^{-1}$'s are in $[-\frac{p-1}{2}, -1]$,
 we have

\hspace{5.25cm} $0 > s \ge -\frac{p(p-1)}{2}$.

\noi Since both $p$ and $q$ are primes, $q \ge p$, and $s = 0$ in
$GF(q)$, we conclude that $q = p$.

Now consider the case that $|I_-| \le p-1$. In this case, we have
$|I_+| \ge p-1$. Applying Lemma (2.1), there is a subset $J$ of
$I_+$, such that

\hspace{5.25cm} $\sum_{i \in J}{x_i}^{-1} = -1$.

\noi holds in $GF(p)$. This implies that $M$ has a
circuit-hyperplane corresponding to $\vec{x}(J)$. Since $A_2$ is
also a special standard representation of $M$, we conclude that
the equality

\hspace{5.25cm} $\sum_{i \in J}{y_i}^{-1} = -1$

\noi holds in $GF(q)$. That is, the equality

\hspace{5.25cm} $\sum_{i \in J}{x_i}^{-1} + 1 = 0$

\noi holds in both $GF(p)$ and $GF(q)$.

Consider the sum $s = \sum_{i \in J}{x_i}^{-1} + 1$ in
$\textbf{Z}$.
 Since all values of $x_i^{-1}$'s are in $[1, \frac{p-1}{2}]$,
 we have

\hspace{5.25cm} $\frac{(p-1)^2}{2}+1 \ge s > 1 $.

\noi Since both $p$ and $q$ are primes, $q \ge p$, and $s = 0$ in
$GF(q)$, we conclude that $q = p$. The main theorem follows
immediately.  \qed

\section{Discussion}

\hspace{.6cm}First we would like to point out that the bound
$2p-1$ is sharp for every prime number $p$. It is easy to prove
the following proposition:

\noi {\bf (4.1) Proposition} {\it Suppose that $M$ is an $n$-spike
representable over $GF(p)$, and $A$ is a special standard
representation of $M$. Let the diagonal of $A$ be $\vec {x} = (
x_1, x_2, \ldots, x_n)$. Suppose that

\hspace{1.2cm}{\rm (1)} $n =2p-2$, and

\hspace{1.2cm}{\rm (2)} $x_1 = x_2 = \ldots = x_p = -1$, and

\hspace{1.2cm}{\rm (3)} $x_{p+1} = x_{p+2} = \ldots = x_{2p-2} =
1$.

\noi Then $M$ is represented by the same matrix $A$ over every
prime field $GF(q)$ with $q \ge p$. }

Characteristic sets of a matroids had been an interesting topic in
matroid theory. The main theorem and proposition (4.1) provide new
and interesting examples for the topic. Readers may also discover
that some typical examples of this topic are in fact spikes.
Besides the Fano and none-Fano matroids, the famous matroids $L_p$
constructed by Lazarson \cite{La58} are also spikes.

An interesting problem related to the main result of this paper
is:

\noi {\bf (4.2) Problem} {\it What is the lower bound L(p) such
that every $GF(p)$-representable $n$-spike with $n < L(p)$ is also
representable over some fields with characteristic other than
$p$?}

We do not have the sharp bound for the above problem at current
time. Our research shows that $L(p)$ is a number between $\lfloor
\log_2(p + 2)+1 \rfloor$, and $\lfloor \log_2(p + 2) \rfloor +
\lfloor \log_2[4(p+2)/3] \rfloor$. However, the argument is
somehow complicated and considered not interesting for our
readers. We instead present the following proposition that is
related to this problem:

\noi {\bf (4.3) Proposition} {\it Suppose that $p$ is an odd prime
number, $M$ is an $n$-spike representable over $GF(p)$, and $A$ is
a special standard representation of $M$. Let the diagonal of $A$
be $\vec {x} = ( x_1, x_2, \ldots, x_n)$. Let $q = \lfloor \log_2p
\rfloor$. Suppose that

\hspace {1.2cm} {\rm (1)} $n =2q+2$,

\hspace {1.2cm} {\rm (2)} $x_1 = -1$,

\hspace {1.2cm} {\rm (3)} $x_{2i}^{-1} = -2^{i-1}$, and
$x_{2i+1}^{-1} = 2^{i-1}$, for $i \in \{ 1, 2, \ldots, q \}$, and

\hspace {1.2cm} {\rm (4)} $x_{2q+2}^{-1} = -2^{q}$.

\noi Then $M$ is only representable over fields of characteristic
$p$. }

\noi \pf Suppose that $F$ is a finite field such that $M$ is
$F$-representable. Suppose that $A^{'}$ is a special standard
representation of $M$ over $F$, and its diagonal is $\vec {y} =
(y_1, y_2, \ldots, y_n)$. By considering the circuit-hyperplane
corresponding to $\{ x_1 \}$, we deduce by applying (2.4) that
$y_1 = -1$. Similarly, we have $y_2 = -1$. Next consider the
circuit-hyperplane corresponding to the vector $( x_1, x_2, x_3
)$. We conclude again by applying (2.4) that $y_3 = x_3 = 1$. Now
switch to consider the circuit-hyperplane corresponding to $(x_3,
x_4 )$. We this time conclude that $y_4^{-1} = x_4^{-1} = -2$. For
$k \ge 3$, by considering circuit-hyperplanes corresponding to
$(x_1, x_{2k-2}, x_{2k-1})$ and $(x_3, x_5, \ldots, x_{2k-1},
x_{2k})$ alternatively, it is not hard to derive that

\hspace {5.cm} $y_i^{-1} = x_i^{-1}$, for each $i \in [1,n]$.

\noi Since $q = \lfloor \log_2p \rfloor$, there is a subset $J$ of
$\{ 2, 4, \ldots, 2q+2 \}$ such that

\hspace {5.2cm}$\sum_{i \in J}{x_i}^{-1} = -p$ in $\textbf{Z}$.

\noi By considering the circuit-hyperplane corresponding to
$\vec{x}(J \bigcup \{ 1 \})$, we conclude that

\hspace {5.2cm}$\sum_{i \in J}{y_i}^{-1} = 0$.

\noi The last equality implies that equation

\hspace {5.2cm} $\sum_{i \in J}{x_i}^{-1} = -p = 0$

\noi holds in both fields $GF(p)$ and $F$. Therefore, $F$ must
have characteristic $p$, and the proposition follows.

\bibliographystyle{amsalpha}

\end{document}